\documentclass[11pt, final]{amsart}
\usepackage{amssymb}

\theoremstyle{plain}
\newtheorem{theorem}{Theorem}
\newtheorem{lemma}{Lemma}

\theoremstyle{remark}

\def\C{\mathbb C}
\def\ker{\operatorname{ker}}
\def\rng{\operatorname{rng}}

\begin{document}
\title[Some Linear Preserver Problems on $B(H)$]
{Some Linear Preserver Problems on $B(H)$ Concerning Rank and Corank}
\author{LAJOS MOLN\'AR}
\address{Institute of Mathematics\\
         Lajos Kossuth University\\
         4010 Debrecen, P.O.Box 12, Hungary}
\email{molnarl@math.klte.hu}
\thanks{  This research was supported from the following sources:\\
          1) Joint Hungarian-Slovene research project supported
          by OMFB in Hungary and the Ministry of Science and
          Technology in Slovenia, Reg. No. SLO-2/96,\\
          2) Hungarian National Foundation for Scientific Research
          (OTKA), Grant No. T--016846 F--019322,\\
          3) A grant from the Ministry of Education, Hungary, Reg.
          No. FKFP 0304/1997}
\begin{abstract}
As a continuation of the work on linear maps between operator algebras
which preserve certain subsets of operators with finite rank, or finite
corank, here we consider the problem
inbetween, that is, we treat the question of preserving operators
with infinite rank and infinite corank. Since, as it turns out,
in this generality our preservers cannot be written in a nice form
what we have got used to when dealing with linear preserver problems,
hence we restrict
our attention to certain important classes of operators like
idempotents, or
projections, or partial isometries. We conclude the paper with a result
on the form of linear maps which preserve the left ideals in $B(H)$.
\end{abstract}
\maketitle

\section{Introduction}

Linear preserver problems represent one of the most active research
topics in matrix theory (see the survey paper \cite{LiTsing}). In the
last
decade considerable attention has been also paid to similar questions
in infinite dimension, that is, to linear preserver problems
on operator algebras (see the
survey paper \cite{BrSe}). In both cases, the problem is to
characterize those linear maps on the algebra in question which leave
invariant a given subset, or relation, or function. One of
the most important such questions concerns the rank. This is because in
many cases preserver problems can be reduced to the problem of
rank preservers. Therefore, it is not surprising that a lot of work has
been done on such preservers (see, for example, \cite{Beas, Loewy}
for the finite dimensional case and \cite{Hou, OmlSem} for the
infinite dimensional case as well as the references therein).
In our recent paper \cite{MolGyorSem}, we considered, among other
things,
the very similar problem of corank preservers which problem deserves
attention, of course, only in the infinite dimensional case. If $H$ is a
(complex) infinite dimensional Hilbert space, denote by $B(H)$ the
algebra of
all bounded linear operators acting on $H$. The result
\cite[Theorem
3]{MolGyorSem} reads as follows. Let $\phi:B(H)\to B(H)$ be a
bijective linear map which is weakly continuous on norm bounded sets.
If $\phi$ preserves the corank-$k$ operators in both directions,
then there exist invertible operators $A,B\in B(H)$ such that
$\phi$ is of the form $\phi(T)=ATB$ $(T\in B(H))$.

Now, it seems to be a natural question to consider the problem of such
preservers which are "inbetween"
rank preservers and corank preservers, that is, to determine those
linear maps
which preserve the operators with infinite rank and infinite corank.
We say that an operator $A\in B(H)$ has infinite rank and infinite
corank if the (Hilbert space) dimensions of $\overline{\rng A}$ and
$\overline{\rng A}^\perp$ are both infinite. Here, $\rng A$ stands for
the range of $A$.
We consider separable Hilbert spaces since in this case
there is only one sort of infinite dimension.
Unfortunately, the preservers above do not have such a nice form
which we have got used to when dealing with linear preserver
problems. Namely, there exist preservers of the above kind
which cannot be expressed in terms of multiplications by fixed operators
and, possibly, by transposition.
To see this, let $\psi:B(H)\to
B(H)$ be a linear map with norm less than 1 whose range consists of
finite rank operators. Then it follows from a basic Banach algebra
fact that the linear map $\phi$ defined by
$$
\phi(T)= T-\psi(T) \qquad (T\in B(H))
$$
is a bijection of $B(H)$ onto itself, and it is easy to check that
$\phi$ preserves
the operators with infinite rank and infinite corank in both directions
(observe that this map preserves the Fredholm index as well which
preserver problem might also seem to be natural after discussing
corank preservers).
So, in order to have one of the desired nice forms for our preservers we
should somehow modify
the problem by, for example, restricting the set of operators with
infinite rank and infinite corank which we want preserve. This is
exactly what we are doing here considering the important
sets of idempotents, projections and partial isometries,
respectively. In the last result of the paper we describe the linear
bijections of $B(H)$ which preserve the left ideals in both directions.
As it will be clear from the proof, this problem is also connected with
the problem of rank preservers.

Let us fix the concepts and notation that we shall use throughout.
By a projection we mean a self-adjoint idempotent in $B(H)$. An element
$W\in B(H)$ is called a partial isometry if it is an isometry
on a closed subspace of $H$ and 0 on its orthogonal
complement. Algebraically, $W$ can be characterized by the
equation $WW^*W=W$.
We say that the operators $A,B\in B(H)$ are orthogonal to each other if
$A^*B=AB^*=0$. This means that the ranges of $A$ and $B$ as well as the
orthogonal complements of their kernels are orthogonal to
each other.
If $x,y\in H$, then $x\otimes y$ denotes the operator defined by
$(x\otimes y)z=\langle z,y\rangle x$ $(z\in H)$.
In what follows $F(H)$ stands for the ideal of all finite rank operators
in $B(H)$.

\vfill

\section{Results}

We begin with the description of all linear bijections $\phi$ of $B(H)$
which
preserve the partial isometries of infinite rank and infinite corank
in both directions (this means that $W$ is a partial isometry with
infinite rank and infinite corank if and only if so is $\phi(A)$).

\begin{theorem}\label{T:partialiso}
Let $H$ be a separable infinite dimensional Hilbert space. Let
$\phi:B(H)\to B(H)$ be a linear bijection which preserves the partial
isometries of infinite rank and infinite corank
in both directions. Then there exist unitary operators $U,V\in B(H)$
such that $\phi$ is either of the form
\[
\phi(T)=UTV \qquad (T\in B(H))
\]
or of the form
\[
\phi(T)=UT^{tr}V \qquad (T\in B(H))
\]
where ${}^{tr}$ denotes the transpose with respect to an arbitrary but
fixed complete orthonormal sequence in $H$.
\end{theorem}

In the proof we shall use the following two auxiliary results.

\begin{lemma}\label{L:partialiso1}
Let $T,S\in B(H)$ be partial isometries with $S=ST^*S$. Then we have
$TT^*S=S$ and $ST^*T=S$.
\end{lemma}

\begin{proof}
Denote $Q=TS^*$.
Since $SS^*$ and $T^*T$ are projections, we compute
\[
SS^*=ST^*SS^*TS^*=Q^*(SS^*)Q\leq
Q^*Q= S(T^*T)S^*\leq SS^*.
\]
This implies $Q^*Q=SS^*$. In particular, we obtain
$\| Q\|\leq 1$ (in fact, the norm of $Q$ is either 0 or 1). But $Q$
is an idempotent. Indeed, we have
\[
Q^2=TS^*TS^*=T(ST^*S)^*=TS^*=Q.
\]
So, $Q$ is a contractive idempotent. It is easy to see that this implies
that $Q$ is a self-adjoint idempotent, that is, a projection. To verify
this, pick arbitrary elements $x \in \ker Q$ and $y\in \rng Q$.
Then we have
\[
\| y\|^2 \leq \| \mu x+y\|^2 \qquad (\mu \in \C).
\]
An elementary argument shows that this implies that $x\perp y$. Hence
the kernel and the range of $Q$ are orthogonal to each other and this
verifies that $Q$ is a projection. Now, from $Q^*Q=SS^*$ we obtain
$Q=SS^*$. Therefore, $TS^*=SS^*$ and, as $SS^*$ is the projection onto
$\rng S$, it follows that the range of $S$
is included in that of $T$. Since $TT^*$ is the projection onto $\rng
T$, we have $TT^*S=S$.
Similarly, from the equality $S^*=S^*TS^*$ one can
deduce $T^*TS^*=S^*$ which is equivalent to $ST^*T=S$.
\end{proof}

\begin{lemma}\label{L:partialiso2}
Suppose that $T,S\in B(H)$ are partial isometries. The
operator $T+\lambda S$ is a partial isometry for every $\lambda \in \C$
with $|\lambda|= 1$ if and only if $T$ and $S$ are orthogonal to each
other.
\end{lemma}

\begin{proof}
Suppose first that
\[
(T+\lambda S) (T+\lambda S)^* (T+\lambda S)=T+\lambda S
\]
holds for every $\lambda\in \C$ with $|\lambda|=1$. Using
the fact that $T, S$ are partial isometries, one can conclude that
\[
0=\lambda^2 S{T}^*S+\lambda(T{T}^*S+S{T}^*T)+\bar \lambda
  T{S}^*T+ S{S}^*T+T{S}^*S.
\]
Since this is valid for every $\lambda\in \C$ of modulus 1, choosing
the particular values $\lambda= 1, -1, i, -i$, it is easy to deduce that
\begin{equation}\label{E:7}
S{T}^*S=0
\end{equation}
\begin{equation}\label{E:8}
T{T}^*S+S{T}^*T=0
\end{equation}
\begin{equation}\label{E:9}
T{S}^*T=0
\end{equation}
\begin{equation}\label{E:10}
S{S}^*T+T{S}^*S=0.
\end{equation}
Multiplying \eqref{E:8} by ${T}^*$ from the left and taking \eqref{E:9}
into account, we obtain ${T}^*S=0$. Similarly, multiplying
\eqref{E:10} by ${S}^*$ from the right and taking \eqref{E:7} into
account, we have $T{S}^*=0$. So, $T$ and $S$ are orthogonal.
As for the reverse implication, if $T,S$ are
mutually orthogonal partial isometries, then it is just a simple
calculation
that $T+\lambda S$ is a partial isometry for every $\lambda \in \C$ of
modulus 1.
\end{proof}

\begin{proof}[Proof of Theorem~\ref{T:partialiso}]
Let $\{ x_1,\ldots, x_k\}\subseteq H$ and $\{ y_1,\ldots, y_k\} \subseteq H$
be two systems of pairwise orthogonal unit vectors. We claim that
the image of the finite rank partial isometry $R=\sum_{j=1}^k x_j
\otimes y_j$ under $\phi$ is also a finite rank partial isometry. Let
$(e_n)$ be an orthonormal sequence in the orthogonal complement of $\{
x_1, \ldots ,x_k\}$ which generates a closed subspace of infinite
codimension. Similarly, let $(f_n)$ be an orthonormal sequence in
$\{ y_1,\ldots, y_k\}^\perp$.
Denote $U=\sum_n e_n \otimes f_n$ and let $V=U+R$. Clearly, $U$ and $V$
are partial isometries of infinite rank and infinite corank.
Moreover, for every $\lambda \in \C$ of modulus 1, the
operator $R+\lambda U=(V-U)+\lambda U$ is also a partial isometry of
infinite rank and infinite corank.
Therefore,
$\phi(V)+(\lambda -1)\phi(U)$ is a partial isometry for every
$\lambda \in \C$ with $|\lambda |=1$. This means that with the notation
$V'=\phi(V), U'=\phi(U)$ we have
\[
(V'+(\lambda-1)U') (V'+(\lambda-1)U')^* (V'+(\lambda-1)U')=
(V'+(\lambda-1)U')
\]
for every $\lambda \in \C$ of modulus 1.
Performing the above operations we obtain a polynomial
in $\lambda, \bar \lambda$ with operator coefficients
which equals 0 on the
perimeter of the unit disc in the complex plane. Just as in the
proof of Lemma~\ref{L:partialiso2}, choosing the particular
values $\lambda=1,-1, i,-i$ we find that the coefficients of the
polynomial in question are all 0. Therefore, we have
\begin{equation}\label{E:1}
-U'+U'{V'}^*U'=0,
\end{equation}
\begin{equation}\label{E:2}
2U'+V'{V'}^*U'+U'{V'}^*V'-U'{U'}^*V'-2U'{V'}^*U'-V'{U'}^*U'=0,
\end{equation}
\begin{equation}\label{E:3}
U'+V'{U'}^*V'-U'{U'}^*V'-V'{U'}^*U'=0,
\end{equation}
and
\begin{multline}\label{E:4}
-2U'-V'{V'}^*U'-V'{U'}^*V'-U'{V'}^*V'+\\
2U'{U'}^*V'+U'{V'}^*U'+2V'{U'}^*U'=0,
\end{multline}
where the left hand sides of \eqref{E:1}, \eqref{E:2}, \eqref{E:3},
\eqref{E:4}
are the coefficients of $\lambda^2$, $\lambda$, $\bar \lambda$ and $1$,
respectively.
From \eqref{E:1} and \eqref{E:2} we deduce
\begin{equation}\label{E:5}
V'{V'}^*U'+U'{V'}^*V'=U'{U'}^*{V'}+{V'}{U'}^*U'.
\end{equation}
We prove that $\phi(R)=V'-U'$ is a partial isometry. Indeed, we compute
\begin{equation}\label{E:6}
(V'-U')(V'-U')^*(V'-U')=
\end{equation}
\[
V'-U'-V'{V'}^*U'-V'{U'}^*V'-U'{V'}^*V'+U'{U'}^*V'+U'{V'}^*U'+V'{U'}^*U'.
\]
From \eqref{E:5} we know that
\[
-V'{V'}^*U'-U'{V'}^*V'+U'{U'}^*V'+V'{U'}^*U'=0.
\]
So, we have to show that
$V'-U'-V'{U'}^*V'+U'{V'}^*U'=V'-U'$.
By \eqref{E:1} we have $U'{V'}^*U'=U'$. It remains to verify that
$V'{U'}^*V'=U'$. From \eqref{E:3} and \eqref{E:5} we infer that
\begin{equation}\label{E:6a}
U'+V'{U'}^*V'-V'{V'}^*U'-U'{V'}^*V'=0.
\end{equation}
By Lemma~\ref{L:partialiso1} it follows that
$V'{V'}^*U'=U'$ and $U'{V'}^*V'=U'$. Now, \eqref{E:6a}
gives $V'{U'}^*V'=U'$. Consequently, the right hand side of the
equation \eqref{E:6} is equal to $V'-U'$ which verifies that
$\phi(R)$ is a partial isometry.

We next prove that $\phi(R)$ has finite rank. By the preserver property
of $\phi$ it is sufficient to prove that $\phi(R)$ has infinite corank.
We have seen that for every $\lambda \in \C$ of modulus 1, the operator
$R+\lambda U$ is a partial isometry of infinite rank and infinite
corank. This implies that
for $R'=\phi(R)$, the operator $R'+\lambda U'$ is a partial isometry
for every $\lambda \in \C$ with $|\lambda |=1$.
According to Lemma~\ref{L:partialiso2}, we obtain that
$R'$ and $U'$ are orthogonal to each other.
Since the range of $U'$ is infinite dimensional, it follows that $R'$ is
of infinite corank which implies that $R'$ is a finite rank partial
isometry.

We next prove that $\phi$ preserves the partial isometries in general.
To see this, let $W$ be
a partial isometry. If it is of finite rank, then there is now nothing
to prove. So, let $W$ be of infinite rank. In that case we have
an orthogonal sequence $(W_n)$ of partial isometries of infinite rank
and infinite corank whose sum is $W$. By the preserver property
of $\phi$ it follows that the operators
$A_n=\phi(W)-\sum_{k=1}^{n+1} \phi(W_k)=\phi(W-\sum_{k=1}^{n+1}W_k)$
and $B_n=\sum_{k=1}^n \phi(W_k)=\phi(\sum_{k=1}^n W_k)$
are partial isometries. Because of the same reason, $A_n +\lambda
B_n$ is a partial isometry for every $\lambda \in \C$ of modulus 1.
By Lemma~\ref{L:partialiso2} this implies that $A_n$ and $B_n$ are
orthogonal to each other.
The statement \cite[Lemma 1.3]{MolZal} tells us that the series of
pairwise orthogonal partial isometries is convergent in the strong
operator topology and its sum is also a partial isometry.
Consider the operators
$A=\phi(W)-\sum_n \phi(W_n)$ and $B=\sum_n \phi(W_n)$.
By the just mentioned result $\sum \phi(W_n)^*$ is strongly convergent
as well, and $\sum_n \phi(W_n)^*=B^*$.
We then also have $\phi(W)^*-\sum_n \phi(W_n)^*=A^*$. Since $A_n$
and $B_n$ are orthogonal for every $n\in \mathbb N$, it is now easy
to verify that $A$ is orthogonal to $B$. The operator $B$ is a partial
isometry. As for $A$,
we know that $(A_n)$ strongly converges to $A$ and, as we have seen,
$(A_n^*)$ strongly converges to $A^*$. It is well-known
that the multiplication is strongly countinuous on the norm-bounded
subsets of $B(H)$. Consequently, we infer that $(A_n A_n^*)$ strongly
converges to $AA^*$ and then that $(A_n A_n^* A_n)$ strongly converges
to
$AA^*A$. Since $A_n$ is a partial isometry for every $n$, we obtain that
$A$ is also a partial isometry. Now, since $\phi(W)$ is the sum of the
mutually orthogonal partial isometries $A$ and $B$, it follows that
$\phi(W)$ is a partial isometry as well. We have assumed that
$\phi^{-1}$
has the same preserver properties as $\phi$. Therefore,
$\phi$ preserves the partial isometries in both directions.
Suppose that $W$ is a maximal partial isometry, that is, suppose that
$W$ is a partial isometry and there is no nonzero partial isometry which
is orthogonal to $W$. If $V\in B(H)$ is a nonzero partial isometry
which is orthogonal to $\phi(W)$, then $V+\lambda \phi(W)$ is a partial
isometry for every $\lambda \in \C$ with $| \lambda |=1$. This gives us
that  $\phi^{-1}(V)+ \lambda W$ is also a partial isometry for every
$\lambda \in \C$ of modulus 1. By Lemma~\ref{L:partialiso2} this results
in the orthogonality of $\phi^{-1}(V)$ and $W$ which is a contradiction.
Consequently, we obtain that $\phi$ preserves the maximal partial
isometries which are precisely the isometries and the coisometries.
It is well-known that the set of all extreme points of the unit ball of
$B(H)$ consists of these operators exactly. So, $\phi$ is a linear
map on $B(H)$ which preserves the extreme points of the unit ball.
The form of all linear maps with this property acting on a von Neumann
factor
was determined in \cite{MolMasc}. The result \cite[Theorem 1]{MolMasc}
says that there is a unitary operator $U\in B(H)$ such that either there
exists a *-homomorphism $\psi:B(H) \to B(H)$ such that
\[
\phi(T)=U\psi(T) \qquad (T\in B(H))
\]
or there exists a *-antihomomorphism $\psi':B(H) \to B(H)$ such that
\[
\phi(T)=U\psi'(T) \qquad (T\in B(H)).
\]
Since our map $\phi$ is bijective, the same must hold for the
corresponding morphism
$\psi$ or $\psi'$ above. Now, referring to folk results on the
form of *-automorphisms and *-antiautomorphisms of $B(H)$, we conclude
the proof.
\end{proof}

We continue with a result of the same spirit on idempotent preservers.

\begin{theorem}\label{T:idempotents}
Let $H$ be a separable infinite dimensional Hilbert space. Suppose that
$\phi:B(H)\to B(H)$ is a linear bijection which preserves the
idempotents of infinite rank and infinite corank in both directions.
Then there is an invertible operator $A\in B(H)$ such that $\phi$ is
either of the form
\[
\phi(T)=ATA^{-1} \qquad (T\in B(H))
\]
or of the form
\[
\phi(T)=AT^{tr}A^{-1} \qquad (T\in B(H)).
\]
\end{theorem}

In the proof we shall use the following lemma which is certainly
well-known and is included here only for the sake of completeness.

\begin{lemma}\label{L:idempotents}
If $P,Q\in B(H)$ are idempotents, then
\begin{itemize}
\item[$(i)$]  $P+Q$ is an idempotent if and only if $PQ=QP=0$;
\item[$(ii)$] $P-Q$ is an idempotent if and only if $PQ=QP=Q$.
\end{itemize}
\end{lemma}

\begin{proof}
It follows from elementary algebraic computations.
\end{proof}

\begin{proof}[Proof of Theorem~\ref{T:idempotents}]
If $P,Q\in B(H)$ are idempotents, then we write $P\leq Q$ if $PQ=QP=P$.
Clearly,
this is equivalent to the condition that $\rng P\subseteq \rng Q$ and
$\ker Q\subseteq \ker P$.
Let us say that an idempotent $P\in B(H)$ is regular if it has infinite
rank and infinite corank.
We prove that for any two regular idempotents $P,Q$ we have $P\leq Q$ if
and only if for every regular idempotent $R\in B(H)$, if $Q+R$ is a
regular idempotent, then so is $P+R$.
The necessity is almost evident. To the sufficiency suppose first that
$\rng P\nsubseteq \rng Q$. Let $x\in H$ be such that $Px=x$ and
$Qx\neq x$. Choose a regular idempotent $R\leq I-Q$ for which $Q+R$ is a
regular idempotent and $Rx\neq 0$ (observe that $(I-Q)x\neq 0$).
Since $P+R$ is an idempotent, we have $PR=RP=0$. It follows that
$0=RPx=Rx$ which is a contradiction. Hence, we have
$\rng P\subseteq \rng Q$. The relation $\ker Q\subseteq \ker P$ can be
proved in a similar manner. Using the above characterization
and the preserver property of $\phi$, we obtain that $\phi$ preserves
the relation $\leq$ between regular idempotents. Now, if $R$ is a finite
rank idempotent, then $R$ can be written in the form $R=Q-P$ with some
regular idempotents $P\leq Q$. Since $\phi(P)\leq \phi(Q)$, it
follows that $\phi(R)=\phi(Q)-\phi(P)$ is also an idempotent. We prove
that $\phi(R)$ is of finite rank.
Choosing a regular idempotent
$P$ with $R\leq P$, it follows that $P-R$ is a regular idempotent and
hence $\phi(P)-\phi(R)$ is also an idempotent. By
Lemma~\ref{L:idempotents} $(ii)$ this implies that
$\phi(R)\leq \phi(P)$. If $\phi(R)$ is not of finite rank, then it is
regular which implies that $R$ is also regular and this is a
contradition. Therefore, using the preserver properties of $\phi$ and
$\phi^{-1}$ we obtain that $\phi$ preserves the finite rank
idempotents in both directions. It is now easy to see that $\phi$ is a
linear bijection of $F(H)$ onto itself. By Lemma~\ref{L:idempotents}
$(i)$, for any idempotents $R,R'\in F(H)$ we have $RR'=R'R=0$ if and
only if $\phi(R)\phi(R')=\phi(R')\phi(R)=0$. Using this property
it is easy to verify that $\phi$
preserves the rank-one idempotents in both directions.
By \cite[Theorem 4.4]{OmlSem} we infer that there is an invertible
bounded
linear operator $A\in B(H)$ such that $\phi$ is either of the form
\[
\phi(T)=ATA^{-1} \qquad (T\in F(H))
\]
or of the form
\[
\phi(T)=AT^{tr}A^{-1} \qquad (T\in F(H)).
\]
Without loss of generality we may assume that $\phi$ is of the first
form and then that $A=I$. We intend to show that $\phi(T)=T$ $(T\in
B(H))$. Let $P\in B(H)$ be a regular idempotent. If $R$ is any finite
rank
idempotent with $R\leq P$, then just as above, we obtain $R=\phi(R)\leq
\phi(P)$. Since $R\leq P$ was arbitrary, it now follows that $P\leq
\phi(P)$.
Since $\phi^{-1}$ has the same preserver property as $\phi$, it follows
that $P\leq \phi^{-1}(P)$. But $\phi$ preserves the order between the
regular idempotents. Hence, we have $\phi(P)\leq P$. Therefore,
$\phi(P)=P$ for every regular idempotent $P$. Since every idempotent of
finite corank is the sum of two regular idempotents, we obtain that
$\phi(P)=P$ holds for every idempotent $P\in B(H)$.
Since every element of $B(H)$ is a finite linear combination of
projections \cite[Theorem 2]{Fil}, we conclude that $\phi(T)=T$
is valid for every $T\in B(H)$. This completes the proof.
\end{proof}

In a similar fashion one can verify the following result
concerning projection preservers.

\begin{theorem}
Let $H$ be a separable infinite dimensional Hilbert space. Suppose that
$\phi:B(H)\to B(H)$ is a linear bijection which preserves the
projections of infinite rank and infinite corank in both directions.
Then there is a unitary operator $U\in B(H)$ such that $\phi$ is
either of the form
\[
\phi(T)=UTU^* \qquad (T\in B(H))
\]
or of the form
\[
\phi(T)=UT^{tr}U^* \qquad (T\in B(H)).
\]
\end{theorem}

Our final result describes the linear bijections $\phi$ of $B(H)$
which preserve the left ideals
in both directions (this means that $\mathcal L\subseteq B(H)$ is a left
ideal if and only if $\phi(\mathcal L)$ is a left ideal).
As it will be clear from the proof,
this problem is also connected with the problem of rank preservers.

\begin{theorem}
Let $H$ be a Hilbert space. Suppose that $\phi:B(H) \to B(H)$ is a
linear bijection preserving the left ideals of $B(H)$ in
both directions. Then there are invertible
operators $A,B\in B(H)$ such that $\phi$ is of the form
\[
\phi(T)=ATB \qquad (T\in B(H)).
\]
\end{theorem}

\begin{proof}
The minimal left ideals of $B(H)$ are precisely the sets
$\{ x\otimes y \, :\, x\in H\}$ for nonzero $y\in H$.
Since $\phi$ clearly
preserves the minimal left ideals of $B(H)$ in both directions, we
easily deduce that
$\phi$ is a linear bijection of $F(H)$ onto itself which preserves the
rank-one operators. By \cite[Theorem 3.3]{OmlSem} (see also \cite{Hou})
it follows that there are linear bijections $A,B:H \to H$ such that $\phi$
is either of the form
\begin{equation}\label{E:20}
\phi(x\otimes y)=Ax \otimes By \qquad (x,y \in H)
\end{equation}
or of the form
\[
\phi(x\otimes y)=Ay \otimes Bx \qquad (x,y \in H).
\]
Since $\phi$ is left ideal preserving,
the second possibility above obviously cannot occur.

We prove that $\phi(I)$ is invertible. First we note the following.
It is true in any algebra with unit that an element fails to have a left
inverse if and only if this element is included in a maximal left ideal
(recall that every proper left ideal is included in a maximal left
ideal).
Therefore, $\phi$ preserves the left invertible elements of $B(H)$ in
both directions. We recall that an operator $S$ in $B(H)$ is left
invertible if and only if $S$ is injective and $S$ has closed range.
Now, let $x,y\in H$ be arbitrary nonzero vectors. Let $\lambda \in \C$.
By Fredholm alternative $x\otimes y-\lambda I$ is injective if and only
if it is surjective. This gives us that
$x\otimes y-\lambda I$ is left invertible if and only if it is
invertible. Since the spectrum of any element in $B(H)$ is nonempty, we
infer that there is a $\lambda \in \C$ for which $x\otimes y-\lambda I$
is not left invertible. Suppose that $x\otimes y$ is not quasinilpotent,
that is, $\langle x,y\rangle \neq 0$. Then the scalar $\lambda$ above
can be chosen to be nonzero.
It follows that $Ax\otimes By-\lambda \phi(I)$
is not left invertible. On the other hand, $\phi(I)$ is left invertible
and hence it is a left Fredholm operator (see \cite[2.3.
Definition, p. 356]{Conway}). But any compact perturbation of a left
Fredholm operator has closed range \cite[2.5. Theorem, p.
356]{Conway}. So, the operator $Ax\otimes By-\lambda \phi(I)$ is not
left invertible but it has closed range. Therefore,
this operator is not injective, that is, there exists a nonzero vector
$z\in H$ such that $\lambda \phi(I)z= \langle z,By \rangle Ax$.
Clearly,
this implies that $Ax\in \rng \phi(I)$. Since $x\in H$ was arbitrary, we
conclude that $H=\rng A\subseteq \phi(I)$ which means that $\phi(I)$ is
surjective. This gives us that $\phi(I)$ is invertible.

We next show that the linear operators $A,B$ in \eqref{E:20} are
bounded.
Let $x,y \in H$. We have seen above that $x\otimes y-\lambda I$ is
not left
invertible if and only if $\lambda \in \sigma (x\otimes y)$, where
$\sigma(.)$ denotes the spectrum. Similarly, by Fredholm alternative
again,
$\phi(I)^{-1}(Ax\otimes By -\lambda \phi(I))$ is not left invertible if
and only if $\lambda \in \sigma (\phi(I)^{-1}Ax \otimes By)$. Since
$\phi$ preserves the left invertible operators in both directions, we
obtain
\[
\sigma (x\otimes y) = \sigma (\phi(I)^{-1}Ax \otimes By).
\]
By the spectral radius formula we have
\[
| \langle x, y\rangle | =| \langle \phi(I)^{-1}Ax, By \rangle|
\qquad (x,y \in H).
\]
Now, an easy application of the closed graph theorem shows that $A,B$
are continuous.

Evidently, we may suppose without any loss of generality
that $A=B=I$.
Let $S\in B(H)$ be invertible and write $C=\phi(S)$. We claim that
$C=S$.
Let $x\in H$ be an arbitrary unit vector. Then $S(I-\lambda x\otimes x)$
has a left inverse if and only if $\lambda\neq 1$. Consequently,
the operator $C-\lambda Sx\otimes x$ is injective for every $\lambda
\neq 1$ and for every unit vector $x\in H$.
Let $z\in H$ be a nonzero vector. Let $y=S^{-1}Cz$ which is also nonzero
since $C=\phi(S)$ is left invertible.
If $\langle z,y\rangle \neq 0$, then choosing $\lambda=\| y\|^2/\langle
z,y\rangle$ we see that
\[
Cz-\lambda (1/\| y\|^2 Sy\otimes y)(z)=Sy-Sy=0.
\]
Since $z$ is nonzero, we deduce $\lambda=1$ which means $\|y\|^2=\langle
z,y\rangle$.
Therefore, for every nonzero vector $z\in H$ we have two possibilities.
Either $\langle z, S^{-1}Cz\rangle=0$ or $\langle
S^{-1}Cz,S^{-1}Cz\rangle=\langle z,S^{-1}Cz\rangle$.
Clearly, the set of all nonzero vectors satisfying the first equality as
well as the set of those ones which satisfy the second equality are both
closed in $H\setminus \{ 0\}$. Since $H\setminus \{ 0\}$ is a connected
set, we infer that either
\[
\langle z, S^{-1}Cz\rangle=0 \qquad (z\in H)
\]
or
\[
\langle S^{-1}Cz,S^{-1}Cz\rangle=\langle z,S^{-1}Cz\rangle \qquad (z\in
H).
\]
The first possibility would imply that $S^{-1}C=0$. Therefore, we have
the
second one which can be reformulated as $(S^{-1}C)^*(S^{-1}C)=S^{-1}C$.
This shows that $S^{-1}C$ is
a projection. But, on the other hand, it is left invertible, and hence
we have $S^{-1}C=I$ which results in $\phi(S)=C=S$.
Since $B(H)$ is linearly generated by the set of all invertible
operators, it follows that $\phi$ is the identity on $B(H)$. This
completes the proof.
\end{proof}

It is easy to see that in the proof above we used only the
preservation of
two extreme kinds of left ideals, namely, that of the minimal
ones and that of the maximal ones.
Preserving minimal left ideals
is connected with the problem of rank preservers. On the other hand,
preserving maximal left ideals
is connected with the problem of left invertibility preservers.
Because of the great
interest in linear maps preserving invertibility in one direction, or
in both directions,
it might be interesting to consider the "one-sided" analogues
of those problems.

\bibliographystyle{amsplain}

\end{document}